\documentclass[conference]{IEEEtran}
\usepackage{graphicx} 
\usepackage{optidef}
\usepackage{breqn}
\usepackage{todonotes}
\usepackage{dirtytalk}
\usepackage{bm}
\usepackage{bbold}
\usepackage{setspace} 
\usepackage{pgfplots}
\pgfplotsset{compat=1.15}
\usepackage{mathrsfs}
\usetikzlibrary{arrows}
\usepackage{amssymb,amsmath,amsthm}
\usepackage{cancel}
\usepackage{tabularx}

{}{}

\newcommand{\cD}{\mathcal{D}}

\newcommand{\cM}{\mathcal{M}}

\newcommand{\cT}{\mathcal{T}}
\newcommand{\cU}{\mathcal{U}}

\newcommand{\cW}{\mathcal{W}}
\newcommand{\cX}{\mathcal{X}}

\newcommand{\E}{\mathbb{E}}

\newcommand{\norm}[1]{\left\lVert #1 \right\rVert}

\theoremstyle{definition}

\definecolor{dtsfsf}{rgb}{0.8274509803921568,0.1843137254901961,0.1843137254901961}
\definecolor{qqqqff}{rgb}{0,0,1}
\definecolor{sexdts}{rgb}{0.1803921568627451,0.49019607843137253,0.19607843137254902}
\definecolor{rvwvcq}{rgb}{0.08235294117647059,0.396078431372549,0.7529411764705882}
\IEEEoverridecommandlockouts
\usepackage{cite}
\usepackage{hyperref}
\usepackage{enumitem}

\usepackage[caption=false,font=footnotesize]{subfig}

\usepackage{algpseudocode}
\usepackage{algorithm}
\usepackage{graphicx}
\usepackage{textcomp}
\usepackage{xcolor}
\def\BibTeX{{\rm B\kern-.05em{\sc i\kern-.025em b}\kern-.08em
    T\kern-.1667em\lower.7ex\hbox{E}\kern-.125emX}}
\usepackage{fancyhdr}
\begin{document}

\title{Stochastic Optimization for Resource Adequacy in Capacity Markets with Storage and Renewables}

\author{\IEEEauthorblockN{Baptiste Rabecq, Andy Sun}
\IEEEauthorblockA{\textit{Operations Research Center,}
\textit{MIT}\\
Cambridge, MA, USA\\
(brabecq, sunx)@mit.edu}
\and
\IEEEauthorblockN{Feng Zhao, Tongxin Zheng, Xiaochu Wang, Yufan Zhang}
\IEEEauthorblockA{\textit{ISO-New England}\\
Holyoke, MA, USA\\
(fzhao, tzheng, xiwang, yuzhang)@iso-ne.com}
}

\maketitle

\begin{abstract} The integration of storage and renewable resources fundamentally alters resource-adequacy analysis. Because storage couples decisions across time, it invalidates the traditional reliability models that are based on time-independent capacity demand curves. Moreover, renewables introduce temporally correlated intermittency. To address this, we formulate the capacity procurement problem as a two-stage stochastic program, where the capacity decision is made in the first stage, while the expected unserved energy is evaluated by a second-stage dispatch problem that considers uncertainties such as generator failures via Markov chains, temporally correlated renewable output, and stochastic load. 
We implement the resulting stochastic capacity procurement (SCP) model on a New England system with 305 generators, including conventional, renewable, and storage units. Using the stochastic decomposition (SD) algorithm, we solve the SCP with up to 20,000 Monte Carlo samples, each representing a six-month trajectory of more than 4,300 hours of uncertainty across all units.
We analyze the convergence behavior of SD and show that convergence for the stochastic program happens faster than reliable estimation of the reliability metrics, which require more samples than are used in typical stochastic programs. These results show that chronologically detailed Monte Carlo sampling can be integrated into capacity procurement optimization in a computationally tractable manner, enabling reliability evaluation with controlled statistical accuracy at realistic system scales.
\end{abstract}


\section{Introduction}

\begin{table}[t]
\small{
\setlength{\tabcolsep}{2pt}
\caption{Notation table for sets, variables, and parameters}
\label{tab:notations}
\centering
\begin{tabular}{|c|p{0.56\columnwidth}|}
\hline
\textbf{Notation} & \textbf{Description} \\
\hline
\multicolumn{2}{|c|}{\textbf{Sets}} \\
\hline
$\cT$ & Set of hourly time steps, of length $T$.\\
$G, R, S$ & Set of conventional generators, renewable generators and energy storage systems\\
$\cD, \cW, \cM$ & Set of days, weeks, months of length $D$, $W$, respectively\\
$\cT_d, \cT_w, \cT_m \subset \cT$ & Set of hours in month $m$, week $w$ and day $d$, respectively.\\
$\cX$ & Set of Capacity Decisions, polyhedral\\
\hline
\multicolumn{2}{|c|}{\textbf{Optimization Variables}} \\
\hline
$x$ & Capacity decision\\
$LS_t$ & Load shed in time period $t$\\ 
$z_t$ & Total power output in time period $t$\\ 
$p_{i, t}$ &Power output of generator $i \in R \cup G$ in time period $t$\\
$e_{s, t}$ & Battery state of charge of storage system $s$ in time period $t$\\
$p^+_{s, t}, p^-_{s, t}, \eta^+_{s}, \eta^-_{s}$ & Charging and discharging power, efficiency of storage system $s$ in time period $t$\\  
\hline
\multicolumn{2}{|c|}{\textbf{Economic parameters}} \\
\hline
$c$& Capacity bids of the generators (\$/kW-month)\\
$VOLL$ & Value of Lost Load (\$/kWh)\\
\hline
\multicolumn{2}{|c|}{\textbf{Physical Parameters}} \\
\hline
$K^{\text{day}}_g$, $K^{\text{week}}_g$, $K^{\text{month}}_g {\in} {[0,1]}$ & Proportion of energy that can be generated by generator $g$ in a day, week or month, respectively. \\
$H_s$ & Duration of storage $s$, in hours.\\
\hline
\end{tabular}
}
\vspace{-2.1em}
\end{table}

In deregulated electricity markets, capacity markets were introduced to achieve long-term resource adequacy \cite{joskow_capacity_2008}.
They treat resource adequacy as a trade-off between the cost of supplying capacity and the value of avoided unserved energy. Generators incur investment and operating costs to access the grid, while consumers suffer an unobserved disutility when load is shed, typically represented by administratively defined downward-sloping demand curves \cite{zhao_constructing_2018}. 
These demand curves are often calibrated to probabilistic reliability criteria such as Loss of Load Expectation (LOLE), or Expected Unserved Energy (EUE) \cite{billinton_reliability_1996}, using the “1-in-10 years” standard. 

With increasing renewable penetration and more rapid variations in demand, assessing the contribution of individual generation units to system reliability becomes increasingly difficult. Recent work from ESIG and DOE argues that, in high-renewable systems, loss-of-load risk is no longer confined to a few peak hours, and that adequacy assessments must evaluate risk across many hours and seasons \cite{doe2024future, stenclik2021redefining}. In response, demand curves are now typically built using Monte Carlo tools that simulate outage, load, and renewable uncertainty \cite{ge_vernova_ra}. However, in systems with significant storage, reliability is inherently intertemporal: a resource’s ability to serve load at a given hour depends on the feasibility of its entire prior state-of-charge trajectory. Scalar capacity metrics and hour-by-hour availability models fail to enforce these long-horizon feasibility constraints and can therefore misrepresent the contribution of storage to system reliability.

This paper formulates resource adequacy (RA) as a stochastic capacity procurement problem (SCP): a two-stage stochastic optimization problem that explicitly accounts for storage dynamics and chronologically detailed uncertainty. First-stage decisions determine generation and storage capacity. In the second stage, adequacy is evaluated through Monte Carlo sampling of long-horizon uncertainty trajectories, where each sampled realization is assessed by solving a load-shedding minimization problem subject to dispatch and state-of-charge constraints. This formulation guarantees that reliability is measured using feasible operating trajectories rather than approximations.

While closely related to the stochastic capacity expansion literature \cite{go_assessing_2016,newlun_adaptive_2021,valencia_zuluaga_parallel_2024}, our setting differs in two critical respects. First, resource adequacy outcomes are driven by temporally correlated rare events that are poorly captured by pre-sampled representative days or small scenario sets, leading to under-representation of loss-of-load conditions. Second, convergence of a finite-scenario stochastic program does not guarantee accurate estimation of reliability metrics such as LOLE or EUE. Indeed, any finite-sample approximation produces a downward-biased estimate of expected unserved energy, and statistical error in reliability criteria decreases only with increasing sample size \cite{shapiro_lectures_2021}. The resulting capacity mix may still produce inaccurate reliability estimates because of insufficient sampling of rare events.

We address this challenge by explicitly separating optimization error from statistical error in our SCP model. We adopt a stochastic decomposition (SD) algorithm \cite{higle_finite_1994} that dynamically samples uncertainty trajectories during the optimization process, allowing capacity decisions to be optimized while assessing the statistical accuracy of reliability metrics.

The paper makes three contributions. First, we propose a storage-aware formulation of SCP that enforces long-horizon intertemporal feasibility. Second, we demonstrate that this formulation can be solved at realistic scales using SD, even with tens of thousands of Monte Carlo scenarios. Third, we provide a statistical assessment of reliability metrics at the optimal solution, showing the distinction between optimization convergence and reliability estimation accuracy. We present numerical results on a stylized ISO New England system with hourly resolution on a six-month horizon and tens of thousands of Monte Carlo scenarios.

Demand curves are an adequate tool to generate market price signals. By contrast, our centralized optimization approach does not directly address this, which is a potential limitation. 
A thorough treatment of market design and pricing implications lies beyond the scope of this paper.

\section{Model}
In this work, we take the central planner's point of view, solving a risk-neutral surplus maximization problem \cite{zhao_constructing_2018, guo_incentivizing_2023}. 
We introduce a compact form for the SCP model:
\begin{align} \label{eq:extensive_problem}
    \min_{x \in \cX} c^\top x + V(x)
\end{align}
where $x$ is the capacity decision, $c^\top x$ is the total capacity payment, and $V(x)$ is the expected cost of unserved energy. Once the capacity decision $x$ is fixed, we evaluate the reliability it delivers to the grid, dependent on uncertainties in system operation. We model three sources of uncertainty:

\begin{enumerate}[leftmargin=*]
    \item \textbf{Generator failures.}
    Following \cite{billinton_reliability_1996}, each unit $i$'s availability is a two-state, time-homogeneous, discrete-time Markov chain on $\{0,1\}$. When $\xi^A_{i, t} = 0$, unit $i$ is unavailable at hour $t$, and when $\xi^A_{i, t} = 1$, it can generate electricity to its full capacity. Reliability is parameterized by the Forced Outage Rate (FOR$_i$) and Mean Time To Repair (MTTR$_i$). Define the per-period repair probability $\lambda_i := 1/\mathrm{MTTR}_i$ and per-period failure probability $\mu_i := \frac{\mathrm{FOR}_i}{(1-\mathrm{FOR}_i)\mathrm{MTTR}_i}$. The availability process $(\xi^A_{i,t})$ has transition matrix
    \begin{align}
        P_i = \begin{pmatrix}
            1-\lambda_i & \lambda_i \\
            \mu_i & 1-\mu_i
        \end{pmatrix}.
    \end{align}
    For each unit $i$, a scenario is a sample path $(\xi^A_{i, 1}, \ldots, \xi^A_{i, T})$ of this Markov chain over the horizon $1, \ldots, T$.
    \item \textbf{Renewable generation.}
    For a renewable unit $r$ and time period $t$, we model generation uncertainty via a capacity factor $\xi^{\mathrm{CF}}_{r,t} \in [0,1]$, which scales the plant’s generation capacity proportionally. To capture temporal correlation within a day, we work with daily capacity-factor profiles. 
    We construct from historical data $K$ representative daily profiles $\{(\hat \xi^{\text{CF}}_{r, t, k})_{t \in 1, \ldots, |\cT_d|}, k \in \{1, \ldots, K\} \}$, and define a discrete distribution $(\pi_{r,1},\ldots,\pi_{r,K})$ over these profiles. 
    For each day $d$ and unit $r$, we draw a profile index $\kappa_{r,d} \in \{1,\ldots,K\}$ independently with $\mathbb{P}(\kappa_{r,d}=k) = \pi_{r,k}$, and set $\xi^{\mathrm{CF}}_{r,t} := \hat\xi^{\mathrm{CF}}_{r,\kappa_{r,d},t}$ for all $t$ in day $d$.
    \item \textbf{Load uncertainty.} 
    We model load uncertainty via a multiplicative factor. Let $(L_t)_{t=1}^T$ be a fixed historical load profile and let $\xi_L \sim \mathrm{Unif}[0.8,1.2]$. The stochastic load process is $\xi^{\mathrm{Load}}_t := \xi_L L_t, \quad t = 1,\ldots,T.$
\end{enumerate}
We combine all randomness into a vector $\xi := (\xi^A,\xi^{\mathrm{CF}},\xi^{\mathrm{Load}})$, and define $(\Xi,\mathcal{F},\mathbb{P})$ a probability space supporting $\xi$. 
We assume that generator failures, renewable profiles, and load are mutually independent, so that the joint law $\mathbb{P}$ factorizes as a product measure of the marginals. 

Following \cite{zhao_constructing_2018}, we assume that the expected cost of unserved load can be written as: 
\begin{align}
     V(x) = VOLL\cdot \E[U(x, \xi)],
\end{align}
where $U(x, \xi)$ is the unserved energy in scenario $\xi$ with capacity decision $x$. For every $x, \xi$, $U(x, \xi)$ can be found as the value function of the following second-stage dispatch problem:
\begingroup
\allowdisplaybreaks
\begin{subequations}\label{eq:second-stage}
\begin{IEEEeqnarray}{rll}
U(x,\xi)&= \min \sum_{t \in \cT} LS_{t} \label{eq:obj}\\
\text{s.t.}&\, LS_t = \xi^{\text{LOAD}}_t - z_t, & \forall t \in \cT \label{eq:balance}\\
&z_{t} {=}\! \begin{aligned}[t]&\sum_{g \in G} p_{g,t} {+}\! \sum_{r \in R} p_{r,t}  \\&\qquad{+}\! \sum_{s \in S} (p_{s,t}^- {-} p_{s,t}^+),\end{aligned} & \forall t \in \cT, \label{eq:total_gen}\\ 
&0 \le p_{g,t} \le \xi_{g,t}^{\mathrm{A}} x_g, & \forall g \in G, t \in \cT \label{eq:gen_limits}\\
&0 \le p_{r,t} \le \xi_{r,t}^{\mathrm{A}} \xi_{r,t}^{\mathrm{CF},k} x_r, & \forall r \in R, t \in \cT \label{eq:ren_limits}\\
&0 \le p_{s,t}^+ \le \xi_{s,t}^{\mathrm{A}} x_s, & \forall s \in S, t \in \cT \label{eq:storage_power_charge}\\
&0 \le p_{s,t}^- \le \xi_{s,t}^{\mathrm{A}} x_s, & \forall s \in S, t \in \cT \label{eq:storage_power_discharge}\\
&0 \le e_{s,t} \le H_s x_s, & \forall s \in S, t \in \cT \label{eq:storage_bounds}\\
&e_{s,1} = 0,& \forall s \in S, \label{eq:storage_init}\\
&e_{s, t+1} {=} e_{s,t} + \eta_s^+ p_{s,t}^+ {-} \frac{p_{s,t}^-}{\eta_s^-}, & \forall s \in S, t\in \cT
\label{eq:storage_dynamics}\\
&\sum_{t \in \mathcal{T}_{d}} p_{g,t} \le K^{\mathrm{day}}_{g} x_g |\cT|_d, & \forall g \in G, d\in \cD \label{eq:fuel_daily}\\
&\sum_{t \in \mathcal{T}_{w}} p_{g,t} \le K^{\mathrm{week}}_{g} x_g |\cT|_w, & \forall g \in G, w\in \cW \label{eq:fuel_weekly}\\
&\sum_{t \in \mathcal{T}_m} p_{g,t} \le K^{\mathrm{month}}_{g} x_g|\cT|_m, & \forall g \in G, m \in \cM\label{eq:fuel_monthly}\\
&LS_{t} \ge 0, & \forall t \in \cT. \label{eq:load_shedding}
\end{IEEEeqnarray}
\end{subequations}
\endgroup

The objective~\eqref{eq:obj} minimizes total load shed. Constraints~\eqref{eq:balance}--\eqref{eq:total_gen} enforce power balance,~\eqref{eq:gen_limits}--\eqref{eq:storage_power_discharge} impose unit and storage power limits,~\eqref{eq:storage_bounds}--\eqref{eq:storage_dynamics} describe storage state-of-charge dynamics, and~\eqref{eq:fuel_daily}--\eqref{eq:fuel_monthly} impose daily, weekly, and monthly energy limits as a proxy for fuel replenishment.

\section{Methods}
Problem~\eqref{eq:second-stage} is a linear program with fixed recourse and the capacity decision $x$ only appears in the right hand side. As a result $U(x, \xi)$ is a convex piecewise linear function of $x$, for almost every $\xi$. Moreover, for any $x \in \cX$ and almost every $\xi$, setting $LS = \xi^{\text{LOAD}}$ and setting all other variables to $0$ yields a feasible solution. Since all variables are bounded above and below, we conclude that for every $x$ and almost every $\xi$, $U(x, \xi) < \infty$. Therefore, problem~\eqref{eq:second-stage} has relatively complete recourse. That is, there exists a (single) dual feasible set $\Pi$, $T$ and $h$ measurable maps such that: 
\begin{align}\label{eq:dual_problem}
    U(x, \xi) = \max_{\theta \in \Pi} \theta^\top (h(\xi) - T(\xi) x),
\end{align}
and there exists a feasible solution for every $x$ and almost every $\xi$.

Problem~\eqref{eq:extensive_problem} is therefore a two-stage stochastic convex program, where the second-stage function $V(x)$ is the expectation of a linear program with relatively complete recourse. While $V(x)$ is not available in closed form, we approximate it by Monte Carlo sampling: given $N_k$ scenarios $(\xi^i)_{i=1, \ldots, N_k}$, we form the sample-average estimator
\begin{align}\label{eq:sample_average}
    \hat V_{N_k}(x) = \frac{1}{N_k} \sum_{i=1}^{N_k}U(x, \xi^i),
\end{align}
which is an unbiased estimator of $V(x)$.
We use the stabilized SD algorithm \cite{higle_finite_1994}, which can be considered a sampling version of a proximally stabilized Stochastic Benders algorithm. Contrary to applying Benders decomposition to a fixed scenario set, SD draws new samples over the course of the algorithm and is designed to approximate the true expectation $V(x)$ rather than a particular finite-sample problem. 

At every iteration $k$, the algorithm constructs a piecewise linear lower approximation $\hat f_k$ of the sample average $\hat V_{N_k}$. Let $\cU_k := \{\alpha_j^k, \beta_j^k\}_{j=1, \ldots, k}$ such that $\hat f_k(x) = \max_{1 \leq j \leq k} \alpha_j^k +  {\beta_j^k}^\top x$. 
Given a proximal penalty $\rho$, and a proximal center $\bar x_k := x_{i_k}$ for some $i_k \leq k$, the iterate $x_{k+1}$ is obtained by solving the following convex Quadratic Program:
\begin{align}\label{eq:epi_model}
x_{k+1} = \arg\min_{x \in \cX} \quad&  c^\top x + \hat f_k(x)+ \frac{\rho}{2} \norm{x - \bar x_k}_2^2.
\end{align}

The details of the cutting plane construction in the SD algorithm are presented in~\ref{alg:sd}. We refer the interested reader to \cite{higle_finite_1994} for a more precise description.

\begin{algorithm}[h]
    \caption{Stabilized Stochastic Decomposition}\label{alg:sd}
    \begin{algorithmic}[1]
    \small
        \State \textbf{Initialize:} $E_0 \gets \emptyset$, $x_0 \in \cX$, $r \in (0, 1)$, $\rho \in \mathbb R^{++}$, $i_{0} = 0$, $\cU_0 = \emptyset$
        \For{$k \in \{1, \ldots, n\}$}, 
            \State Sample $\xi^{N_{k-1} +1}, \ldots, \xi^{N_{k}}$ from $\mathbb P$.
            \State \textbf{Solve dispatch problems:} Solve \eqref{eq:second-stage} at $x_{k}$ and $x_{i_{k-1}}$ for every sample, record optimal duals $\Theta_{k}:= \{\theta^{k}_{N_k +1}, \ldots, \theta^{k}_{N_{k +1}}, \bar \theta^{k}_{N_k +1}, \ldots, \bar \theta^k_{N_{k +1}}\}$
            \State $\Pi_{k} \gets \Pi_{k-1} \cup \Theta_{k}$.
            \State \textbf{Update Models:}   
            \For{$j \in \{1, \ldots, k\}$}
                \If{$j \in \{i_k, k\}$}
                 \State $\theta_t^{j} {:=} {\arg}{\max}_{\Pi_{k}} \theta^\top  (h(\xi^t) {-} T(\xi^t) x_{j})$, $\forall t \in [N_{k+1}]$
                 \State $\alpha_{j}^{k} {\gets} \frac{1}{N_{k}} \sum_{t=1}^{N_{k}} {\theta^k_t}^\top h(\xi^t)$, $\beta_{j}^{k}{\gets} \frac{-1}{N_{k}} \sum_{t=1}^{N_{k}} {\theta^j_t}^\top T(\xi^t)$
                \Else
                \State $\alpha_{j}^{k} \gets \frac{k-1}{k} \alpha_j^{k-1}$, $\beta_{j}^{k} \gets \frac{k-1}{k} \beta_j^{k-1}$ 
                \EndIf            
            \EndFor
            \State $\mathcal U_k \gets (\alpha_j^k, \beta_j^k)_{j \leq k}$
            \If {$
    \hat f_{k}(x_{k}) - \hat f_{k}(\bar x_{k-1}) \leq r (\hat f_{k-1}(x_{k}) - \hat f_{k-1}(\bar x_{k-1}))$}{ $i_{k} \gets k$}
            \Else { $i_{k} \gets i_{k-1}$}
            \EndIf
            \State Get $x_{k+1}$ from \eqref{eq:epi_model}
        \EndFor
    \end{algorithmic}
\end{algorithm}

We emphasize that the cuts are explicitly recomputed only at the current proximal center and the current iterate; at all other points, the corresponding cuts are simply decayed. Because problem \eqref{eq:second-stage} has a nonnegative optimal value, this decay preserves the property that the model is a statistical lower approximation of the true objective $V(x)$. 
The described dual caching procedure allows for increasing the accuracy of the model compared to batch first-order methods (all samples are used to generate cuts), while keeping the iteration cost low compared to Benders decomposition. 

The quantity $\Delta_k:= \hat f_k({\bar x_k}) - \hat f_k(x_{k+1})$ is the model gap, representing the gap between the current minimizer of the model and the proximal center. A small gap is a good indicator of good in-sample convergence.      
\begin{figure*}[t]
    \centering
    \subfloat[Model gap $\Delta_k$]{
        \includegraphics[width=0.3\textwidth]{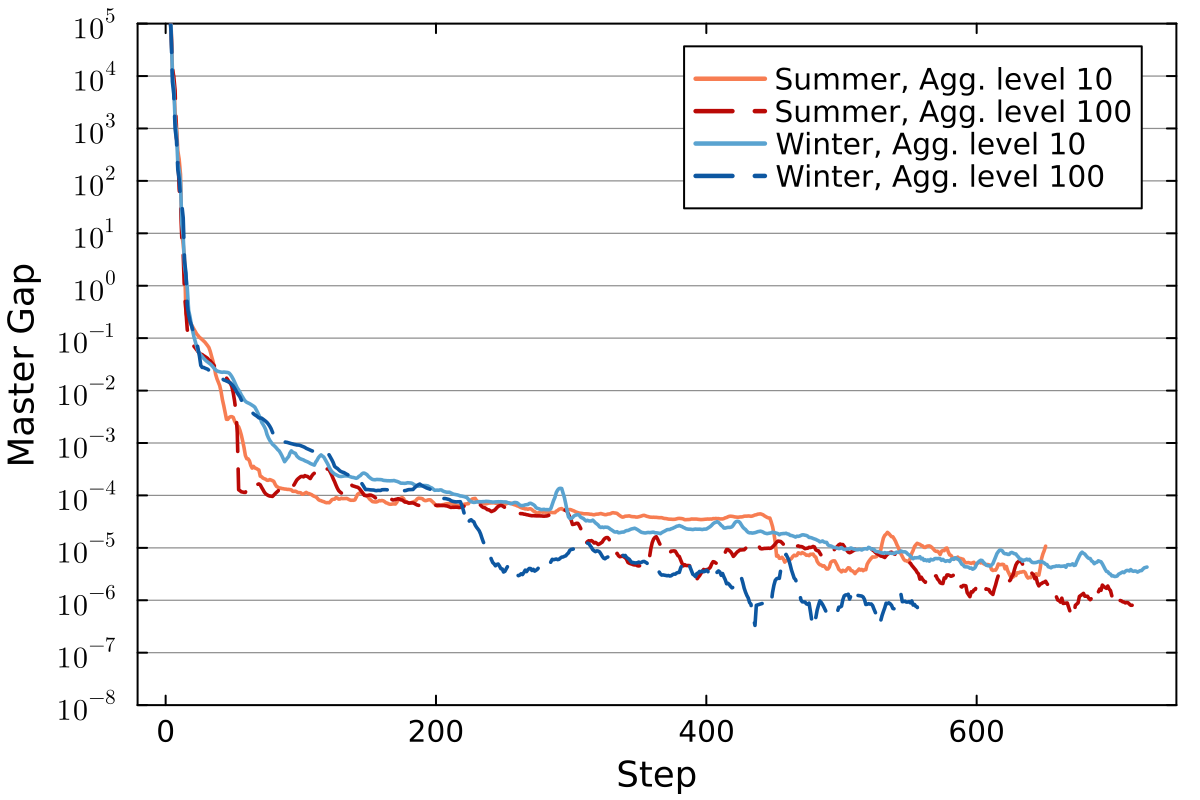}%
        \label{fig:convergence-gap}
    }
    \hfil
    \subfloat[Incumbent objective value]{
        \includegraphics[width=0.3\textwidth]{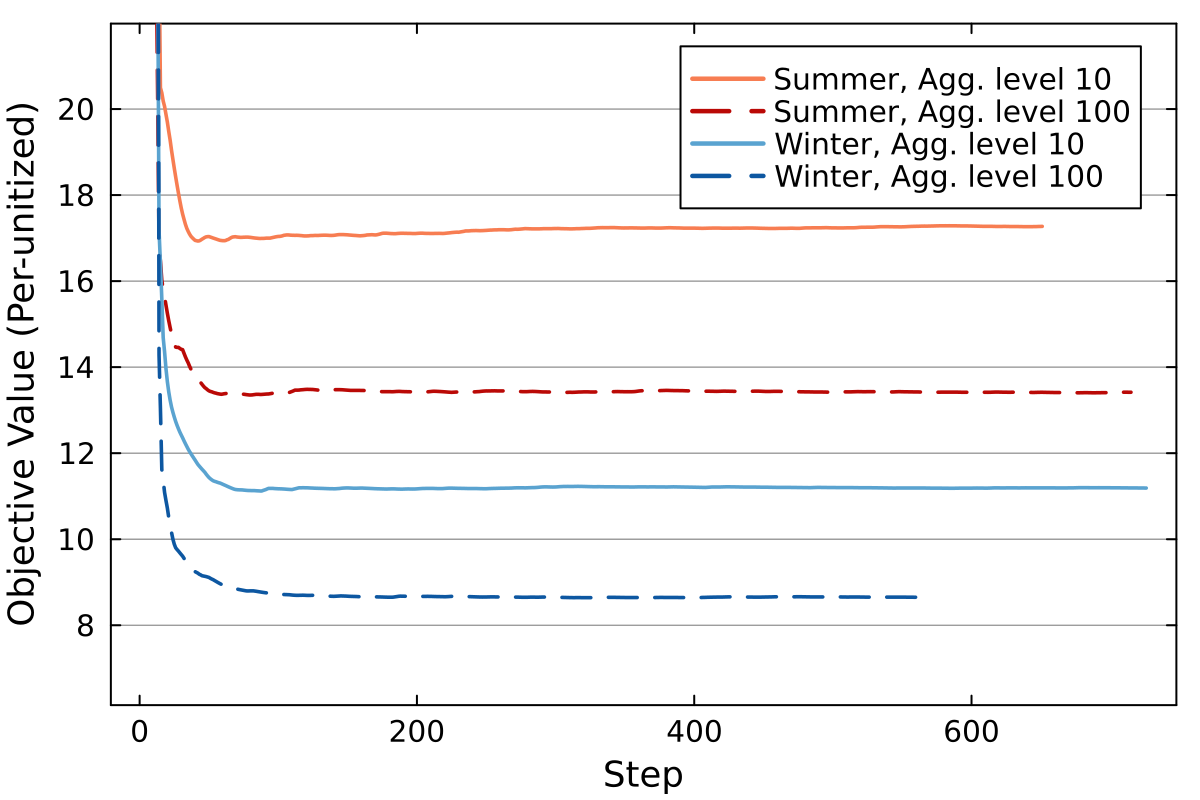}%
        \label{fig:convergence-obj}
    }
    \hfil
    \subfloat[Expected cost of unserved energy]{
        \includegraphics[width=0.3\textwidth]{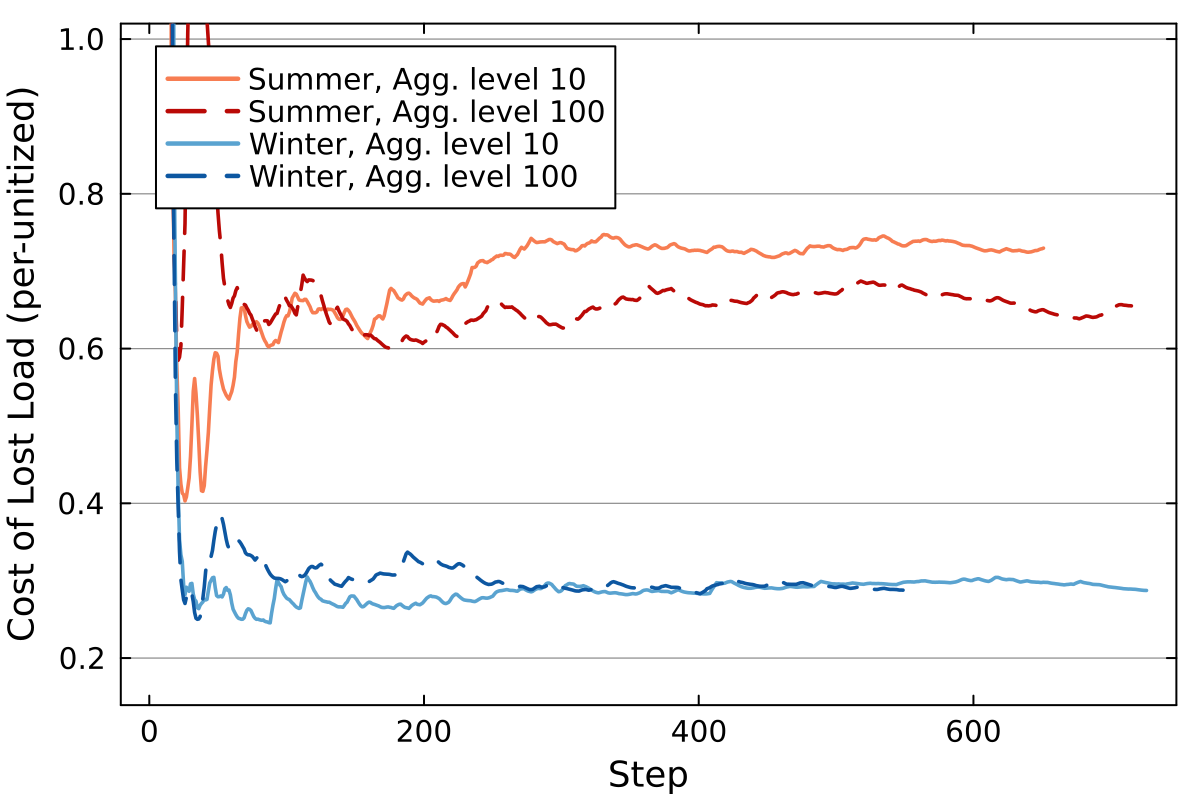}%
        \label{fig:convergence-eue-cost}
    }

    \caption{Convergence diagnostics of the SD algorithm: (a) model gap $\Delta_k$,
    (b) incumbent objective value, and (c) expected cost of unserved energy, all
    as functions of the number of iterations for the summer and winter capacity auction, at aggregation level $10$ ($5$ hours) MW and $100$ MW ($2$ hours).}
    \label{fig:convergence-behavior}
    \vspace{-1em}
\end{figure*}
Evaluating the expected cost of lost load $V(x)$ is a tail estimation problem. Under acceptable reliability standards, $U(x, \xi) =0$, in almost all scenarios, so only a small fraction of realizations are informative about $V$.
Coupled with the high-dimensional uncertainty, this makes finite-sample estimators prone to overfitting in low to medium sample regimes, and achieving a small generalization error typically requires large scenario sets. 

\section{Experiments}

We solve the stochastic capacity procurement model using data from the ISO New England (ISO-NE) system.
Feasible capacity decisions are based on cleared capacity from the ISO-NE 18th Forward Capacity Auction~\cite{ISO_NE_FCM_Reports}, with bidding prices derived from EIA cost data~\cite{EIA_AEO2023_CostPerf} and assumed capacity market revenue shares. Forced–outage uses ISO-NE's Equivalent Forced Outage Rate demand statistics by technology \cite{isone_nerc_gads_eford_2024}, and mean time to repair values from NREL’s Regional Energy Deployment System~2.0 documentation~\cite{NREL_ReEDS2_2025}.

Renewable generation profiles are obtained from one year of hourly solar and wind data from the NREL National Solar Radiation Data Base~\cite{sengupta_national_2018} and the NREL Bias-Corrected High-Resolution Rapid Refresh dataset~\cite{buster_bias_2024}, respectively. We convert these series into capacity factor profiles using Direct Normal Irradiance for solar and the \texttt{windpowerlib} Python package for wind \cite{haas_windpowerlib_2024}. We capture intra-day temporal correlation using $k$-medoids clustering to obtain representative daily profiles and probabilities, sampling five profiles per season.

We clear capacity in two distinct auctions: one summer auction,  from May~1 to October~30, and one winter auction, from November~1 to April~30. We solve the problems with hourly dispatch decisions on the corresponding time horizons of approximately $4300$ hours. 
To limit the size of the problem, we aggregate units of the same technology type when their capacity is below a prescribed threshold, taken as $100$~MW and $10$~MW. With a $100$~MW threshold, the system comprises $107$ units in total, including 10 energy storage units, $2$ solar units, $4$ wind units, and $91$ conventional units.
With a 10~MW threshold, the system comprises $305$ units, including $17$ energy storage units, $10$ solar units, $9$ wind units, and $269$ conventional units. 
We use a Value of Lost Load of \$$100,000$/MWh. According to market practice, we set a market cap based on the net Cost of New Entry, set to \$$14.5$/kW-month.

All experiments are conducted on the Engaging MIT computing cluster, using AMD EPYC 9654 96-core processors. We set the SD batch size to $b = 32$ scenarios per iteration, solved in parallel at every iteration of the algorithm. The choice of $b$ balances per-iteration computational speed with efficiency: larger batches slow each iteration, but smaller batches require more iterations to reach comparable solution quality.

\section{Results}
Throughout the runs, we monitor the in-sample standard deviation of our $EUE$ estimator and stop when it is on the order of $25\%$ of the estimated $EUE$. On our system, this number is reached when approximately $20,000$ samples are used. We therefore use this total number of samples as a stopping criterion for our experiments.  
For the $100$MW aggregated system, this happens after 2 hours, and $5$ hours for the $10$MW aggregated system.

\subsection{Convergence of the algorithm}
Figure~\ref{fig:convergence-behavior} reports convergence diagnostics for the SD algorithm at aggregation level $10$ and $100$ MW over $5$ and $2$ hours of wall-clock time, respectively. Panel~\ref{fig:convergence-gap} shows the model gap $\Delta_k$ as a function of the iteration index. For both aggregation levels, the gap decreases rapidly and falls below $10^{-4}$ after a few hundred iterations, indicating that the master problem is well solved relative to the current sample information. 
Panel~\ref{fig:convergence-obj} displays the incumbent objective value $\hat f_k(\hat x_k)$; it decreases substantially during the first $\sim 100$ iterations, after which improvements are marginal and subsequent iterations produce small fluctuations around a stable level. The difference in objective value between aggregation levels is due to a different per-unit scaling.  Panel~\ref{fig:convergence-eue-cost} shows the corresponding estimated expected cost of unserved energy, which also stabilizes after the first $300$ iterations, and is the same for both aggregation levels. We observe that the lower-dimensional runs ($100$ MW aggregation) converge faster than the higher-dimensional runs. 

A small model gap and visually stable objective values indicate satisfactory convergence of the SD procedure, but they do not certify optimality for Problem~\eqref{eq:extensive_problem}. Because the estimates $\hat f_k(\hat x_k)$ are computed from finite samples, the incumbent decisions $x_k$ continue to move slightly even when $\Delta_k$ is very small, reflecting residual sampling noise. In this work, we assess the final solution $x^*$ via OOS statistical validation. 

\subsection{Statistical validation of the solution}
We evaluate the final solution $x^*$ returned by the algorithm both in-sample (IS) (i.e. during the algorithm execution) and out-of-sample (OOS), using the $10$MW-thresholded system. Table~\ref{tab:oos-validation} reports the estimated reliability metrics under the IS scenario set used during optimization, together with OOS estimates obtained from an independent Monte Carlo simulation of $20,000$ samples. Across all metrics, the IS and OOS confidence intervals overlap, and the IS–OOS differences are of the same order as these intervals, with no consistent sign. This provides no statistically significant evidence of systematic optimization-induced bias at the current Monte Carlo resolution.
The relative CI widths in Table~\ref{tab:oos-validation} range between $0.15$ and $0.30$, with larger values for EUE than for LOLE. Combined with the small EUE levels (below 100~MWh/year), this confirms that expected costs are driven by rare, high-impact events rather than by typical operating conditions.

Finally, Table~\ref{tab:oos-validation} shows that at the optimal solution, LOLE levels are relatively conservative for this stylized system. Indeed, the LOLE is around $45\%$ of the standard value of $0.1$ days per year.

\begin{table}[t]
\centering
\caption{IS vs OOS performance of $x^*$. LOLE, EUE with 95\% CI}
\label{tab:oos-validation}
\begin{tabular}{@{}lccc@{}}
\hline
Season/Metric & IS mean & OOS mean & \shortstack{Rel. OOS\\CI width} \\
\hline
Summer EUE   & $91.70 \pm 11.14$ & $82.24 \pm 9.11$ & 0.22 \\
Summer LOLE  & $0.045 \pm 0.003$ & $0.057  \pm 0.004$ & 0.15 \\
Winter EUE   & $36.82 \pm 5.48$  & $44.34  \pm 6.78$  & 0.30 \\
Winter LOLE  & $0.046 \pm 0.004$ & 0.049  $\pm 0.004$ & 0.16 \\
\hline
\end{tabular}
\vspace{-1em}
\end{table}
\begin{figure}[t]
    \centering
    \subfloat[Summer\label{fig:capacity_mix_summer}]{
        \includegraphics[width=0.472\linewidth]{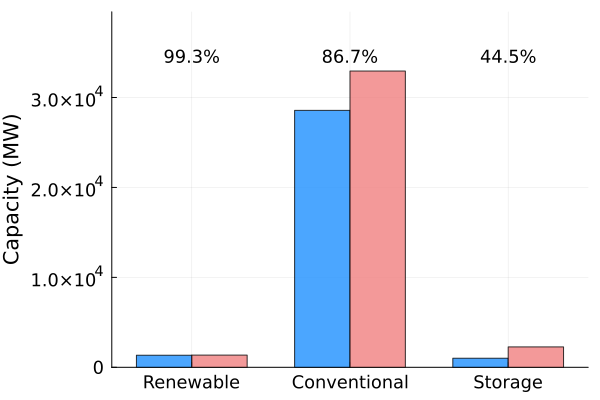}
    }\hfill
    \subfloat[Winter\label{fig:capacity_mix_winter}]{
        \includegraphics[width=0.472\linewidth]{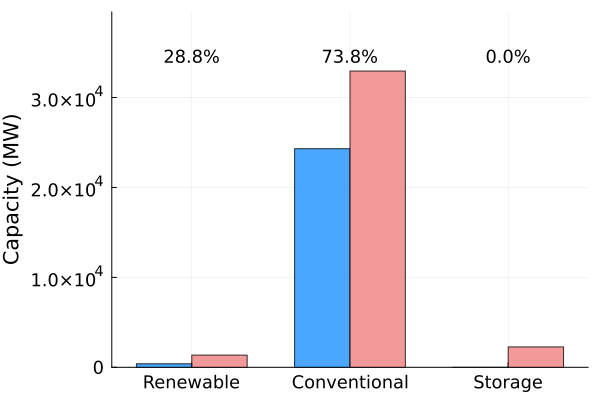}
    }
    \caption{Comparison between summer and winter capacity mix of $x^*$ (MW). The blue and pink bars represent the cleared capacity and total available capacity, respectively.}
    \label{fig:capacity_mix}
    \vspace{-2em}
\end{figure}

\begin{figure}[t]
    \centering
    \subfloat[Summer\label{fig:ls_distribution_summer}]{
        \includegraphics[width=0.472\linewidth]{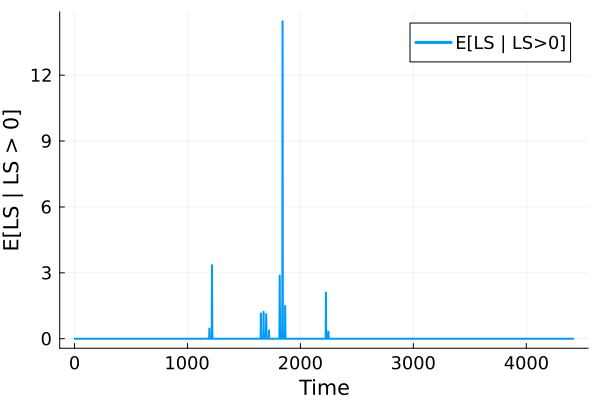}
    }\hfill
    \subfloat[Winter\label{fig:ls_distribution_winter}]{
        \includegraphics[width=0.472\linewidth]{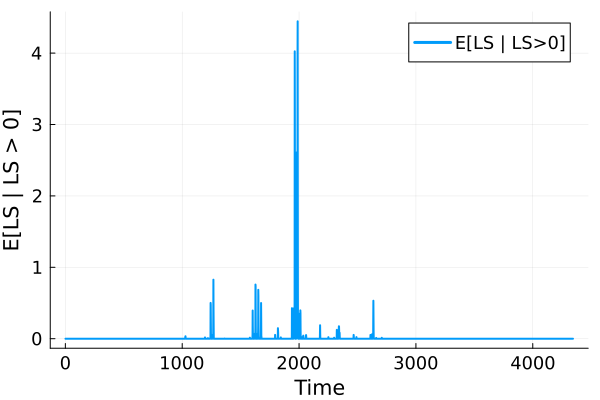}
    }
    \caption{Load Shedding Distribution for the summer and winter auction (MWh)}
    \label{fig:ls_distribution}
\end{figure}

\subsection{Solution and seasonal mix}
Figure~\ref{fig:capacity_mix} displays the capacity mix $x^*$ in Summer and Winter. The generators are divided into three categories: ``Renewable", denoting wind and solar generators, ``Energy Storage", and ``conventional", denoting the remaining ones. We observe that more capacity is contracted in the summer auction than in the winter auction. This is expected, as the capacity is primarily driven by the peak load, which is higher in the summer. In the winter auction, there is little renewable generation and no storage. In the summer, the renewable capacity allocation increases, and some storage capacity is allocated. 
In Figure~\ref{fig:ls_distribution}, we observe that Loss of Load is concentrated on a relatively small set of peak hours, as expected. In some scenarios, the generator's fuel constraints are binding, showing dependence of the EUE on long-term coupling constraints.

\section{Conclusion}
In this work, we propose a two-stage stochastic program for achieving resource adequacy in capacity markets with significant storage and renewables. We solve the resulting large-scale stochastic program using stochastic decomposition with chronological Monte Carlo sampling. Numerical experiments on an ISO-NE system show that the proposed method is computationally efficient with statistically reliable performance. 

\bibliography{main} 

@book{billinton_reliability_1996,
    title = {Reliability {Evaluation} of {Power} {Systems}},
    author = {Billinton, Roy and Allan, Ronald N.},
    year = {1996},
    publisher = {Springer US},
}

@article{higle_finite_1994,
    title = {Finite master programs in regularized stochastic decomposition},
    volume = {67},
    doi = {10.1007/BF01582219},
    language = {en},
    number = {1},
    journal = {Mathematical Programming},
    author = {Higle, Julia L. and Sen, Suvrajeet},
    year = {1994},
    pages = {143--168},
}

@article{sengupta_national_2018,
    title = {The {National} {Solar} {Radiation} {Data} {Base} ({NSRDB})},
    volume = {89},
    issn = {1364-0321},
    doi = {10.1016/j.rser.2018.03.003},
    journal = {Renewable and Sustainable Energy Reviews},
    author = {Sengupta, Manajit and Xie, Yu and Lopez, Anthony and Habte, Aron and Maclaurin, Galen and Shelby, James},
    year = {2018},
    pages = {51--60},
}

@misc{buster_bias_2024,
    title = {Bias {Correcting} {NOAA}'s {High}-{Resolution} {Rapid} {Refresh} ({HRRR}) {Wind} {Resource} {Data} for {Grid} {Integration} {Applications}},
    doi = {10.2172/2479268},
    language = {American English},
    author = {Buster, Grant and Pinchuk, Pavlo and Lavin, Luke and Benton, Brandon and Bodini, Nicola},
    year = {2024},
}

@misc{NREL_ReEDS2_2025,
  author       = {National Renewable Energy Laboratory},
  title        = {Model Documentation --- {ReEDS} 2.0},
  year         = {2025},
  howpublished = {\url{https://nrel.github.io/ReEDS-2.0/model_documentation.html}},
}

@misc{ISO_NE_FCM_Reports,
  author       = {{ISO New England}},
  title        = {Forward Capacity Market Reports --- ISO Express},
  howpublished = {\url{https://www.iso-ne.com/isoexpress/web/reports/auctions/-/tree/forward-capacity-mkt}},
}

@misc{haas_windpowerlib_2024,
  author       = {Haas, Sabine and Krien, Uwe and Schachler, Birgit and Stickler Bot and
                  Zeli, Velibor and Maurer, Florian and Shivam, Kumar and Witte, Francesco and
                  Rasti, Sasan Jacob and Seth and Bosch, Stephen},
  title        = {wind-python/windpowerlib: Update release (v0.2.2)},
  year         = {2024},
  doi          = {10.5281/zenodo.10685057},
}

@techreport{EIA_AEO2023_CostPerf,
  author      = {{U.S. EIA}},
  institution      = {{U.S. EIA}},
  title       = {Cost and Performance Characteristics of New Generating Technologies, {Annual Energy Outlook} 2023},
  year        = {2023},
}

@article{joskow_capacity_2008,
    series = {Capacity {Mechanisms} in {Imperfect} {Electricity} {Markets}},
    title = {Capacity payments in imperfect electricity markets: {Need} and design},
    volume = {16},
    doi = {10.1016/j.jup.2007.10.003},
    number = {3},
    journal = {Utilities Policy},
    author = {Joskow, Paul L.},
    year = {2008},
    pages = {159--170},
}

@article{zhao_constructing_2018,
    title = {Constructing {Demand} {Curves} in {Forward} {Capacity} {Market}},
    volume = {33},
    doi = {10.1109/TPWRS.2017.2686785},
    number = {1},
    journal = {IEEE Transactions on Power Systems},
    author = {Zhao, Feng and Zheng, Tongxin and Litvinov, Eugene},
    year = {2018},
    pages = {525--535},
}

@techreport{stenclik2021redefining,
  author      = {Stenclik, Derek and Bloom, Aaron and Cole, Wesley and
                 Figueroa Acevedo, Armando and Stephen, Gord and Tuohy, Aidan},
  title       = {Redefining Resource Adequacy for Modern Power Systems},
  year        = {2021},
  doi         = {10.2172/1961567},
  institution    = {ESIG},
}

@techreport{doe2024future,
  institution      = {{U.S. DOE}},
  author      = {{DOE}},
  title       = {The Future of Resource Adequacy},
  year        = {2024},
  url         = {https://www.energy.gov/sites/default/files/2024-04/2024%20The%20Future%20of%20Resource%20Adequacy%20Report.pdf},

}

@article{go_assessing_2016,
    title = {Assessing the economic value of co-optimized grid-scale energy storage investments in supporting high renewable portfolio standards},
    volume = {183},
    doi = {10.1016/j.apenergy.2016.08.134},
    journal = {Applied Energy},
    author = {Go, Roderick S. and Munoz, Francisco D. and Watson, Jean-Paul},
    year = {2016},
    pages = {902--913},
}

@inproceedings{newlun_adaptive_2021,
	title = {Adaptive {Expansion} {Planning} {Framework} for {MISO} {Transmission} {Planning} {Process}},
	doi = {10.1109/KPEC51835.2021.9446221},
	urldate = {2025-11-17},
	booktitle = {2021 {IEEE} {Kansas} {Power} and {Energy} {Conference} ({KPEC})},
	author = {Newlun, Cody J. and McCalley, James D. and Amitava, Rajaz and Ardakani, Ali Jahanbani and Venkatraman, Abhinav and Figueroa - Acevedo, Armando L.},
	year = {2021},
	pages = {1--6},
}

@article{valencia_zuluaga_parallel_2024,
    title = {Parallel computing for power system climate resiliency: {Solving} a large-scale stochastic capacity expansion problem with mpi-sppy},
    volume = {235},
    doi = {10.1016/j.epsr.2024.110720},
    journal = {Electric Power Systems Research},
    author = {Valencia Zuluaga, Tomas and Musselman, Amelia and Watson, Jean-Paul and Oren, Shmuel S.},
    year = {2024},
    pages = {110720},
}

@misc{ge_vernova_ra,
    author      = {{GE}},
    title = {Resource {Adequacy} {Software} {\textbar} {GE} {Vernova}},
}

@techreport{isone_nerc_gads_eford_2024,
  author      = {{ISO New England}},
  institution = {{ISONE}},
  title       = {NERC {GADS} {EFORd} Class Averages as used by {ISO} New England},
  year        = {2024},
  type        = {Technical report},
 }

@misc{guo_incentivizing_2023,
    title = {Incentivizing {Investment} and {Reliability}: {A} {Study} on {Electricity} {Capacity} {Markets}},
    language = {en-US},
    author = {Guo, Cheng and Kroer, Christian and Dvorkin, Yury and Bienstock, Daniel},
    year = {2023},
}

@incollection{shapiro_lectures_2021,
author = {Alexander Shapiro},
booktitle = {Lectures on Stochastic Programming: Modeling and Theory, Third Edition},
title = {Chapter 5: Statistical Inference},
pages = {151-221},
doi = {10.1137/1.9781611976595.ch5},
eprint = {https://epubs.siam.org/doi/pdf/10.1137/1.9781611976595.ch5}
}
\end{document}